\documentclass[onecolumn,showpacs,preprintnumbers,amsmath,amssymb]{revtex4-2}

\usepackage{graphicx}
\usepackage{amsmath, amssymb, amsfonts, color}
\usepackage{enumerate, tabularx}
\usepackage{hyperref, url}
\usepackage[margin=1in]{geometry}
\usepackage{caption,subcaption,float}
\usepackage[capitalise]{cleveref}
\newcommand{\edit}[1]{ #1}

\newcommand{\pd}[2]    { \frac{\partial #1} {\partial #2} }

\newcommand{\pdi}[2] { {\partial_#2} #1 }

\newcommand{\abs}[1]{\left| #1 \right|}
\newcommand{\CC}{\mathbb{C}}
\newcommand{\RR}{\mathbb{R}}

\newcommand{\bvec}[1]{\mathbf{#1}}

\newcommand{\uhat}{\hat{u}}
\newcommand{\vard}[2]{\frac{\delta #1}{\delta #2}}
\newcommand{\Ham}{\mathcal{H}}
\newcommand{\En}{\mathcal{E}}
\newcommand{\Mo}{\mathcal{M}}
\newcommand{\nmodes}{K}

\newcommand{\proj}{\mathcal{P}_{\nmodes}}
\newcommand{\Gibbs}{\mathcal{G}}
\newcommand{\bprime}{\beta'}
\newcommand{\sphere}{\mathbb{S}^{2\nmodes - 1}}
\newcommand{\xc}{\xi}
\newcommand{\Sig}{\bvec{\Sigma}}
\newcommand{\diag}{\mathop{\mathrm{diag}}}

\newcommand{\xv}{\bvec{x}}

\newcommand{\Xv}{\bvec{X}}
\newcommand{\Xhv}{\bvec{\hat{X}}}

\newcommand{\xhat}{\hat{x}}

\newcommand{\xhv}{\bvec{\xhat}}

\newcommand{\Ufun}{\bvec{\hat{U}}}

\newcommand{\Const}{Const}

\begin{document}

\title{In search of rogue waves: a novel proposal distribution for parallelized rejection sampling of the truncated KdV Gibbs measure}


\author{Nicholas J.~Moore}
\email{nickmoore83@gmail.com}
\author{Brendan Foerster}
\email{brendancfoerster@gmail.com}
\affiliation{Department of Mathematics, Colgate University, Hamilton, NY 13346, USA}

\begin{abstract}
The Gibbs ensemble of the truncated KdV (TKdV) equation has been shown to accurately describe the anomalous wave statistics observed in laboratory experiments, in particular the emergence of extreme events. Here, we introduce a novel proposal distribution that facilitates efficient rejection sampling of the TKdV Gibbs measure. Within parameter regimes accessible to laboratory experiments and capable of producing extreme events, the proposal distribution generates 1--6 orders of magnitude more accepted samples than does a naive, uniform distribution. When equipped with the new proposal distribution, a simple rejection algorithm enjoys key advantages over a Markov chain Monte Carlo algorithm, include better parallelization properties and generation of uncorrelated samples.
\end{abstract}
\maketitle

\section{Introduction}

Abnormally large water waves --- variously known as anomalous waves, extreme waves, or rogue waves --- have in recent decades captured the interest of the scientific community, ocean practitioners, and the public at large \cite{Haver2004, Muller2005, Garrett2009, Hadj2014, Dem2018, Dem2019, Bolles2019, dudley2019rogue}. Although several mechanisms for their existence have been proposed, one line of research has demonstrated that abrupt changes in bottom topography can trigger anomalous wave activity \cite{trulsen2012laboratory, Viotti2014, Bolles2019, Majda2019, MajdaQi2019, Machine2019, Herterich2019, Trulsen2020, MooreJNS2020, Zheng2020, qi2022anomalous, li2021surface1, li2021surface2, li2021rogue, lawrence2021statistical, mendes2021non, afzal2021propagation}. In particular, laboratory measurements of topography-induced anomalous wave statistics \cite{Bolles2019,MooreJNS2020} were successfully reconciled with a theoretical framework based on statistical and dynamical analysis of the variable-depth Korteweg-De \edit{Vries} (KdV) equation \cite{Majda2019, MajdaQi2019, MooreJNS2020}, a partial different equation (PDE) that describes nonlinear dispersion of shallow waves \cite{Lax1975,Johnson1997,Whitham2011}. Specifically, the theory relies on a Galerkin truncation of KdV, known as the truncated KdV (TKdV) system \cite{MooreJNS2020}. 

The framework exploits the Hamiltonian structure of TKdV to define a Gibbs measure from which wave events can be sampled without the need to directly simulate the underlying PDEs \cite{MooreJNS2020}. 
\edit{In particular, the variable-depth TKdV equation gives rise to separate Gibbs ensembles upstream and downstream of the depth change. The inverse temperatures of these two systems are linked through a statistical matching condition enforced at the depth change \cite{Majda2019,MooreJNS2020}. This framework has been shown to successfully predict outgoing wave statistics as a function of the incoming wave-field characteristics and the size of the depth change.}

\edit{
While the framework has proven successful in predicting system behavior, the task of sampling the Gibbs ensembles presents some challenges.}
Previous studies have employed brute-force sampling by taking the \edit{spectrally} uniform measure as a proposal distribution in a rejection or sampling-importance resampling (SIR) algorithm \cite{MooreJNS2020, SunMoore2023}, or they have employed more sophisticated methods such as Markov chain Monte Carlo (MCMC) \cite{Majda2019,MooreJNS2020}. However, each of these approaches has drawbacks.
\edit{Taking the proposal distribution of a rejection (or SIR) algorithm to be uniform in spectral space} results in an extremely low acceptance rate (or extremely low sample weights) in the parameter regime of greatest interest. Although the MCMC algorithm enjoys higher efficiency, disadvantages include: (1) a number of tunable parameters to calibrate; (2) long ``burn-in'' times and consequently poor parallelization properties; (3) correlation between subsequent samples. An ideal sampling algorithm would combine the easy parallelization and uncorrelated nature of rejection sampling with the high efficiency of MCMC.

Of course, it is possible to improve rejection sampling considerably with a better choice of the proposal distribution. An ideal proposal distribution satisfies two criteria: (1) it approximates the target distribution well, and (2) it is easy to sample. Although the \edit{spectrally} uniform measure satisfies the second requirement, it does not approximate the TKdV Gibbs measure particularly well. On the other hand, Sun \& Moore (2023) \cite{SunMoore2023}, in their efforts to prove theoretical results regarding TKdV wave statistics, constructed a distribution that satisfies both criteria. The goal of the present paper is to employ the distribution of Sun \& Moore (2023) as a proposal distribution in a highly efficient and parallelizable sampling algorithm. 
\edit{We employ this proposal distribution in a {\em rejection} sampling algorithm, although the same distribution could just as easily be employed in an SIR method. To demonstrate the desirable numerical features, it suffices to apply the algorithm to the downstream Gibbs ensemble, assuming that its inverse temperature has already been set by the statistical matching condition \cite{SunMoore2023}. We will therefore consider the constant-depth TKdV system throughout this paper. The application of the algorithm to piecewise-constant depth topography is straightforward and will be done in future work.}

The outline of the paper is as follows. \Cref{Sec:Background} discusses the relevant background information, including the TKdV measure, its Hamiltonian structure, and the canonical-microcanonical Gibbs measure. In \cref{Sec:ProposalDist}, we introduce the proposal distribution, including a recap of the relevant analysis of Sun \& Moore (2023) \cite{SunMoore2023}. \Cref{Sec:Algorithms} outlines the new sampling algorithm.
In \cref{Sec:Results} we present results from numerical tests, including evaluation of the algorithm's performance and physical interpretation of sampled wave fields. \Cref{Sec:Discussion} offers some concluding remarks.

\section{Mathematical preliminaries} 
\label{Sec:Background}

This section provides the mathematical background for the new sampling algorithm. Below, we introduce the KdV and TKdV systems, as well as the associated Gibbs measure. The KdV system derives from the more primitive Euler equations in the long-wave limit, and it retains the physical effects of nonlinearity and dispersion \cite{Lax1975,Johnson1997,Whitham2011}. Analysis of the KdV system complements studies that directly simulate the Euler equations \cite{Viotti2014, ambrose2022numerical}, many of which employ conformal mapping of multiply-connected domains \cite{crowdy2020solving, baddoo2023generalization} to represent obstacles or topographical changes.

\subsection{The Korteweg-De \edit{Vries} equation}

The KdV equation is given by
\begin{align}
\label{KdV}
&u_t + C_3 \, u u_\xc + C_2 \, u_{\xc \xc \xc} = 0
\qquad \text{for } \xc \in [-\pi,\pi] \, ,
\end{align}
where $u(\xc,t)$ is the surface displacement measured in a reference frame traveling with the leading-order wave speed. Boundary conditions are periodic on the normalized domain $\xc \in [-\pi,\pi]$, where $\xc$ is the horizontal location in the traveling frame. All variables $u, \xc, t$ are assumed dimensionless already, and the dimensionless coefficients $C_3$ and $C_2$ characterize the strength of nonlinearity and dispersion respectively \cite{MooreJNS2020}. \edit{Depth variations correspond to variation of the coefficients $C_3$ and $C_2$ with respect to time (since the reference frame is traveling) \cite{Majda2019,MooreJNS2020}. In particular, an abrupt depth change results in distinct Gibbs ensembles upstream and downstream of the depth change, where the downstream inverse temperature is determined by the statistical matching condition \cite{Majda2019,MooreJNS2020}. To demonstrate the features of the new sampling algorithm, it suffices to consider the downstream Gibbs ensemble only, supposing its inverse temperature has already been set by the statistical matching condition. We therefore consider constant values of $C_3$ and $C_2$  throughout this paper. }

The KdV equation possesses a Hamiltonian structure given by
\begin{align}
\label{H3H2}
& \Ham_3[u] = \frac{1}{6} \int_{-\pi}^{\pi} u^3 \, d\xc	\, , \qquad
\Ham_2[u] = \frac{1}{2} \int_{-\pi}^{\pi} \left( \pd{u}{\xc} \right)^2 \, d\xc	\, , \\
\label{KdVHam}
& \Ham[u] = C_2 \, \Ham_2[u] - C_3 \, \Ham_3[u] \, .
\end{align}
\Cref{KdV} can then be written as
\begin{align}
\label{HamStruct}
\pd{u}{t} = \pd{}{\xc} \vard{\Ham}{u} \, ,
\end{align}
where $\pdi{}{\xc}$ is a symplectic operator. Hence, \cref{KdV} is a Hamiltonian system and, consequently, the Hamiltonian, \cref{KdVHam}, is conserved during evolution.

In addition to the Hamiltonian, momentum and energy are conserved under KdV dynamics,
\begin{align}
\label{MomEn}
\Mo[u] \equiv \int_{-\pi}^{\pi} u \, d\xc \, = 0 , \qquad
\En[u] \equiv \frac{1}{2} \int_{-\pi}^{\pi} u^2 \, d\xc = \En_0 \, .
\end{align}
The momentum vanishes, as indicated above, because $u$ is measured as displacement from equilibrium. Meanwhile, the energy has been normalized to the value $\En_0$.
We remark that setting $\En_0 = \pi/2$ recovers the setup of Sun \& Moore (2023) \cite{SunMoore2023}. In the present paper, we will keep $\En_0$ arbitrary for slightly greater flexibility.

\subsection{The truncated KdV system}

We next perform a finite Galerkin truncation of \cref{KdV} \cite{Majda2019, MooreJNS2020, SunMoore2023}. To this end, consider a spatial Fourier representation of the state variable 
\begin{align}
&u(\xc,t) = \sum_{k=-\infty}^{\infty} \uhat_k(t) \, e^{i k \xc} 
=  \sum_{k=1}^{\infty} a_k(t) \cos(k \xc) + b_k(t) \sin(k \xc) \, , \\
\label{uhat}
&\uhat_k = \frac{1}{2} (a_k - i b_k)= \frac{1}{2 \pi} \int_{-\pi}^{\pi} u(\xc,t) \, e^{-i k \xc} \, d\xc \, .
\end{align}
For convenience, we have recorded both the real and the complex Fourier representations, $\uhat_k \in \CC$ and $a_k, b_k \in \RR$. Note that $\uhat_{-k} = \uhat_{k}^*$ since $u(\xc,t)$ is real valued and $\uhat_0 = 0$ due to momentum vanishing.
Next, consider the Galerkin truncation at wavenumber $\nmodes$
\begin{align}
\label{uL}
u_{\nmodes}(\xc,t) = \proj [u] = \sum_{k=-\nmodes}^{\nmodes} \uhat_k \, e^{i k \xc} 
= \sum_{k=1}^{\nmodes} a_k \cos(k \xc) + b_k \sin(k \xc) \, ,
\end{align}
Here, the truncation operator $\proj$ projects an infinite-dimensional function space onto a finite-dimensional one. Note that the formula for the Fourier coefficients \cref{uhat} still holds, even after truncation. Inserting the projected variable, $u_{\nmodes}$, into the KdV equation and applying the projection operator, $\proj$, again where necessary produces the truncated KdV equation (TKdV) \cite{Bajars2013, Majda2019, MooreJNS2020}
\begin{align}
\label{TKdV}
&\pd{u_{\nmodes}}{t} + \frac{1}{2} C_3 \, \pd{}{\xc} \proj \left[ (u_{\nmodes})^2 \right]
+ C_2 \, \frac{\partial^3 u_{\nmodes}}{\partial \xc^3} = 0
\qquad \text{for } \xc \in [-\pi,\pi] \, .
\end{align}
\Cref{TKdV} represents a {\em finite} dimensional dynamical system. The quadratic nonlinearity, $\pdi{}{\xc} \proj \left[ (u_{\nmodes})^2 \right]$, mixes the modes during evolution. The presence of the additional projection operator in this term removes the aliased modes of wavenumber larger than $\nmodes$. Typical values of the cutoff wavenumber used in the previous studies are $\nmodes = $ 8--32 \cite{Majda2019, MooreJNS2020, SunMoore2023, SunMooreBao2023}.

The TKdV equation enjoys nearly the same Hamiltonian structure as KdV, with the only modification being the inclusion of the projection operator,
\begin{align}
\label{TKdVHam}
&\Ham_{\nmodes} = C_2 \, \Ham_2[u_{\nmodes}] - C_3 \, \Ham_3[u_{\nmodes}] \, , \\
&\pd{}{t} {u_{\nmodes}} = \pd{}{\xc} \proj \left[ \vard{\Ham_{\nmodes}}{u_{\nmodes}} \right] \, ,
\end{align}
where now $\pdi{}{\xc} \proj$ is the symplectic operator of interest. 

The system's microstate can either by described in physical space $u_{\nmodes}(\xc, t)$, in complex spectral space, $(\uhat_1, \uhat_2, \cdots, \uhat_{\nmodes}) \in \CC^{\nmodes}$, or in real spectral space $(a_1, a_2, \cdots, a_{\nmodes}, b_1, b_2, \cdots, b_{\nmodes}) \in \RR^{2 \nmodes}$. All are equivalent through \cref{uhat,uL}.
The momentum and energy defined in \cref{MomEn} are also conserved in the truncated system and have the same normalized values $\Mo[u_{\nmodes}] = 0$ and $\En[u_{\nmodes}] = \En_0$. Parseval's identity implies
\begin{equation}
\label{Econd}
\En[u_{\nmodes}] = 2 \pi \sum_{k=1}^{\nmodes} \abs{\uhat_k}^2 = \frac{\pi}{2}
\sum_{k=1}^{\nmodes} a_k^2 + b_k^2 = \En_0 \, .
\end{equation}

\subsection{The Gibbs ensemble}

Following previous studies \cite{Bajars2013, Majda2019, MooreJNS2020, SunMoore2023}, we employ a {\em mixed canonical-microcanonical} Gibbs measure of TKdV; that is, a measure that is canonical in the Hamiltonian and microcanonical in the energy, formally expressed as
\begin{align}
\label{Gibbs0}
d \Gibbs \propto \exp(-\beta \Ham_{\nmodes}) \delta(\En - \En_0) \, .
\end{align}
Above, $\beta$ is the inverse temperature. The exponential dependence with respect to the Hamiltonian is the well-known canonical distribution, which, under suitable conditions, maximizes entropy \cite{MajdaWang2006}. The Dirac-delta term $\delta(\En - \En_0)$ confines the distribution to the compact set $\En = \En_0$, thereby avoiding the far-field, sign-indefinite divergence of $\Ham_3$ and rendering the distribution normalizable  \cite{Abramov2003,Bajars2013,Majda2019,MooreJNS2020,SunMoore2023}. 

This mixed Gibbs measure can be defined rigorously via integration against test functions on the unit hypersphere $\sphere = \lbrace \xv \in \RR^{2\nmodes} \colon \abs{\xv} = 1 \rbrace$. Throughout, we let $\phi: \sphere \to \RR$ represent an arbitrary, measurable test function. First, consider the standard normal distribution on $\RR^{2\nmodes}$,
\begin{equation}
\label{stand_norm}
d \gamma(\xv) :=  \prod_{k=1}^{2\nmodes} \frac{1}{\sqrt{2 \pi}} e^{ {-x_k^2}/{2} } \, dx_k \, .
\end{equation}
The {\em uniform} measure $d\mu_0$ on $\sphere$ can be defined by projecting $d\gamma(\xv)$ onto $\sphere$; that is, for any measureable $\phi: \sphere \to \RR$, let
\begin{equation}
\label{mu0}
\int_{\sphere} \phi(\xhv) \, d\mu_0(\xhv) = \int_{\RR^{2\nmodes}} \phi (\xhv) \, d\gamma(\xv) \, ,
\end{equation}
where $\xhv = {\xv}/{\abs{\xv}} \in \sphere$ represents a unit-vector projection. 
Above, we have chosen a free constant to render $d \mu_0$ a probability density function (p.d.f.). That is, $\int_{\sphere} d\mu_0 = 1$, as can be seen by taking $\phi \equiv 1$ and recalling that $d\gamma$ is a p.d.f.
\edit{We remark that $d\mu_0$ is a uniform measure in {\em spectral} space (not, for example, a uniform distribution of surface displacements).} 

Next, for each unit vector, $\xhv = (\xhat_1,\xhat_2,\cdots,\xhat_{2\nmodes}) \in \sphere$, we can identify a corresponding microstate
\begin{equation}
\label{uhxh}
\uhat_k = \hat{U}_k(\xhv) 
= \sqrt{ \frac{\En_0}{2\pi} } \, (\xhat_k - i \xhat_{\nmodes+k})
\quad \text{for } k=1,2,\cdots,\nmodes
\end{equation}
with the conjugate relation, $\uhat_{-k} = \uhat_{k}^*$, defining the negative modes. By definition, \cref{uhxh} satisfies the energy constraint \cref{Econd}.

With this nomenclature established, we can precisely define the Gibbs measure $d\Gibbs$ from \cref{Gibbs0}. In particular, for any measurable $\phi: \sphere \to \RR$, let
\begin{align}
\label{Gibbs1}
\int_{\sphere} \phi(\xhv) \, d\Gibbs(\xhv) 
&= Z^{-1} \int_{\sphere} \phi(\xhv) \, \exp( -\beta \, \Ham_{\nmodes}[\Ufun(\xhv)]) \, d\mu_0 (\xhv) \\
\label{Gibbs2}
&= Z^{-1} \int_{\RR^{2 \nmodes}} \phi (\xhv) \, \exp( -\beta \, \Ham_{\nmodes}[\Ufun(\xhv)]) \, d\gamma (\xv)
\end{align}
where $Z$ is the partition function (a normalization constant chosen to satisfy the law of total probability) and $\Ufun = (\hat{U}_1, \hat{U}_2, \cdots, \hat{U}_\nmodes) \colon \sphere \to \CC^{\nmodes}$ as given by \cref{uhxh}.
The Gibbs measure is associated with a probability density supported on the hypersphere, $f:\sphere \to \RR$,
\begin{equation}
\label{f_compute}
f(\xhv) = Z^{-1} \, \exp\left( -\beta \, \Ham_{\nmodes}[ \Ufun( \xhv ) ] \right)
\qquad \text{for } \xhv \in \sphere \, .
\end{equation}

We remark that if the inverse temperature vanishes, $\beta = 0$, then the Gibbs measure, $d\Gibbs$ from \cref{Gibbs1}, reduces to the uniform measure, $d \mu_0$, on $\sphere$. Earlier studies have shown that if $\beta \ne 0$, $d\Gibbs$ can differ substantially from the uniform measure and produce anomalous  surface-displacement statistics \cite{MooreJNS2020,SunMoore2023}. More specifically, positive inverse temperature, $\beta >0$, corresponds to the anomalous wave behavior observed in laboratory experiments \cite{Bolles2019,MooreJNS2020}. Our main goal is to devise an algorithm to sample the Gibbs measure \cref{Gibbs1} in this physically relevant and nontrivial regime of $\beta > 0$.

\section{The proposal distribution}
\label{Sec:ProposalDist}

In this section, we introduce the proposal distribution that will enable efficient rejection sampling of the Gibbs measure, \cref{Gibbs1}.
This distribution was first constructed by Sun \& Moore (2023) \cite{SunMoore2023} who used it to establish rigorous results on the conditions that lead to Gaussian versus non-Gaussian surface displacement statistics. The construction is based on analysis of the {\em linear} TKdV system, i.e.~$C_3 = 0$, even though the resulting distribution has relevance to the fully nonlinear system. Here, we briefly recap the derivation of the distribution, with emphasis placed on the intuition behind its construction since the rigor has already been established by Sun \& Moore (2023) \cite{SunMoore2023}.

\subsection{Linear TKdV and change of measure}

Setting $C_3 = 0$ eliminates the only nonlinearity in TKdV, \cref{TKdV}, rendering it a {\em linear} dynamical system. In this case, the Hamiltonian, \cref{TKdVHam}, reduces to its quadratic component $\Ham_{\nmodes} = C_2 \Ham_2$, which, by Parseval's identity, can be written as
\begin{equation}
\label{H2Pars}
\Ham_2 = 2 \pi \sum_{k=1}^{\nmodes} k^2 \abs{\uhat_k}^2 
= \En_0 \sum_{k=1}^{\nmodes} k^2 (\xhat_k^2 + \xhat_{\nmodes+k}^2)
\, .
\end{equation}
To simplify notation, consider a normalized inverse temperature
\begin{equation}
\label{beta_prime}
\bprime =  \En_0 C_2 \nmodes^2 \beta \, ,
\end{equation}
Using \cref{uhxh,Gibbs2,H2Pars,beta_prime}, we can rewrite the main part of the linear TKdV Gibbs measure in terms of the integration variable $\xv \in \RR^{2\nmodes}$,
\begin{equation}
\label{linearGibbs1}
\exp \left( -\beta C_2 \Ham_2 \right) d\gamma(\xv)
= \exp \left( -\frac{\bprime}{\nmodes^2 \abs{\xv}^2} \sum_{k=1}^{\nmodes} k^2 (x_k^2 + x_{\nmodes+k}^2 ) \right) d\gamma(\xv) \, .
\end{equation}

The analysis of Sun \& Moore (2023) \cite{SunMoore2023} builds upon a few observations of this formula. First, for a random variable $\Xv \in \RR^{2\nmodes}$ selected from the standard Guassian distribution $d\gamma$, the law of large numbers (LLN) implies $\abs{\Xv}^2/(2 \nmodes) \to 1$ as $\nmodes \to \infty$, where the convergence is in probability. This convergence might suggest the crude substitution $\abs{\xv}^2 \approx 2 \nmodes$ in \cref{linearGibbs1}, which can be made even more general by introducing an extra degree of freedom $\abs{\xv}^2 \approx 2 \nmodes/\alpha$, where $\alpha > 0$. As it turns out, by selecting $\alpha$ judiciously, this second substitution leads to an approximation of the Gibbs measure that can be justified rigorously \cite{SunMoore2023}.

A brief overview of the derivation is as follows. First, notice that the formula
\begin{equation}
\label{xsq_rec}
\frac{1}{\abs{\xv}^2} = \frac{\alpha}{2 \nmodes} 
+ \frac{1}{2 \nmodes} \left( \frac{2 \nmodes}{\abs{\xv}^2} - \alpha \right)
\end{equation}
holds for any $\alpha$. To offer some intuition, our desired substitution $\abs{\xv}^2 \approx 2 \nmodes/\alpha$ suggests the main part of $1/\abs{\xv}^2$ will be captured by the first term, $\alpha / (2\nmodes)$, and that the second term above will be small.

Next, consider an anisotropic Gaussian measure $\gamma_{\Sig}$ with covariance matrix
\begin{align}
\label{Sig}
&\Sig = \diag \left(\sigma_1^2, \sigma_2^2, \cdots, \sigma_{\nmodes}^2, \sigma_1^2, \sigma_2^2, \cdots, \sigma_{\nmodes}^2 \right) \, , \\
\label{sigma}
&\sigma_k^2 = \frac{1}{1 + \alpha \bprime k^2/\nmodes^3}
\qquad \text{for } k=1,2,\cdots, \nmodes \, .
\end{align}
More explicitly
\begin{align}
\label{dgSig1}
d \gamma_{\Sig}(\xv) 
& = \prod_{k=1}^{\nmodes} \frac{1}{2 \pi \sigma_k^2} \exp \left( -\frac{x_k^2 + x_{\nmodes+k}^2} {2 \sigma_k^2}\right) dx_k dx_{\nmodes+k} \\
\label{dgSig2}
& = C \exp \left(- \frac{\alpha \bprime}{2 \nmodes^3} \sum_{k=1}^{\nmodes} k^2 (x_k^2 + x_{\nmodes+k}^2) \right) \, d\gamma(\xv) \, .
\end{align}
where $C$ is a constant.
In particular, this anisotropic measure $d \gamma_{\Sig}$ is chosen to equal the right-hand-side of \cref{linearGibbs1} when the crude substitution $\abs{\xv}^2 \approx 2 \nmodes/\alpha$ is made. If instead the exact \cref{xsq_rec} is used, then \cref{linearGibbs1} can be rewritten {\em with no approximation made} as 
\begin{equation}
\label{linearGibbs2}
\exp \left( -\beta C_2 \Ham_2 \right) \, d\gamma(\xv)
= C^{-1} \exp \left( - \frac{ \bprime}{ 2\nmodes^3 } \left( \frac{2 \nmodes}{\abs{\xv}^2} - \alpha \right) \sum_{k=1}^{\nmodes} k^2 \left(x_k^2 + x_{\nmodes+k}^2 \right) \right) \, d\gamma_{\Sig}(\xv) \, .
\end{equation}
If $( 2K/\abs{\xv}^2 - \alpha ) \to 0$, as expected from \cref{xsq_rec}, then the right-hand-side of \cref{linearGibbs2} converges to  a constant multiple of $d\gamma_{\Sig}$. This is the main intuition for how the distribution $d\gamma_{\Sig}$ approximates the Gibbs measures in the case of linear TKdV.

Finally, Sun \& Moore (2023) \cite{SunMoore2023} showed that if $\alpha = \alpha^*$ is selected as the root of the function
\begin{equation}
\label{F_alpha}
F(\alpha) := 1 - \frac{\alpha}{\nmodes} \sum_{k=1}^{\nmodes} \frac{1}{1 + \alpha \bprime k^2 / \nmodes^3}
\end{equation}
then for any continuous $\phi: \sphere \to \RR$, the following convergence holds for linear TKdV
\begin{equation}
\label{conv_result}
\int_{\sphere} \phi(\xhv) \, d\Gibbs(\xhv) - \int_{\RR^{2 \nmodes}} \phi(\xhv) \, d\gamma_{\Sig}(\xv) \to 0
\qquad \text{as } \nmodes \to \infty \, ,
\end{equation}
for any fixed $\bprime \ge 0$.
{\em In this precise sense, the anisotropic Gaussian $d \gamma_{\Sig}$, with $\alpha^*$ chosen as the root of \cref{F_alpha}, approximates the linear-TKDV Gibbs measure as the cutoff wavenumber grows large.}

We remark that $\alpha^*$ is not necessarily equal to one, as might have been guessed from the earlier intuitive argument involving the LLN (i.e.~if $\Xv \in \RR^{2\nmodes}$ is selected from the standard Gaussian distribution $d\gamma$, then $\abs{\Xv}^2/(2 \nmodes) \to 1$ as $\nmodes \to \infty$ in probability).
The reason is that, with change of measure \cref{linearGibbs2}, the random variable $\Xv \in \RR^{2\nmodes}$ should now be regarded as coming from the anisotropic Gaussian $d\gamma_{\Sig}$. The value $\alpha = \alpha^*$ is selected self-consistently, so that when the random variable $\Xv$ is sampled from $d\gamma_{\Sig}$ (which depends on $\alpha$), the convergence $2\nmodes/\abs{\Xv}^2 \to \alpha^*$ holds as $\nmodes \to \infty$.

\begin{figure}
\centering
\includegraphics[width=\linewidth]{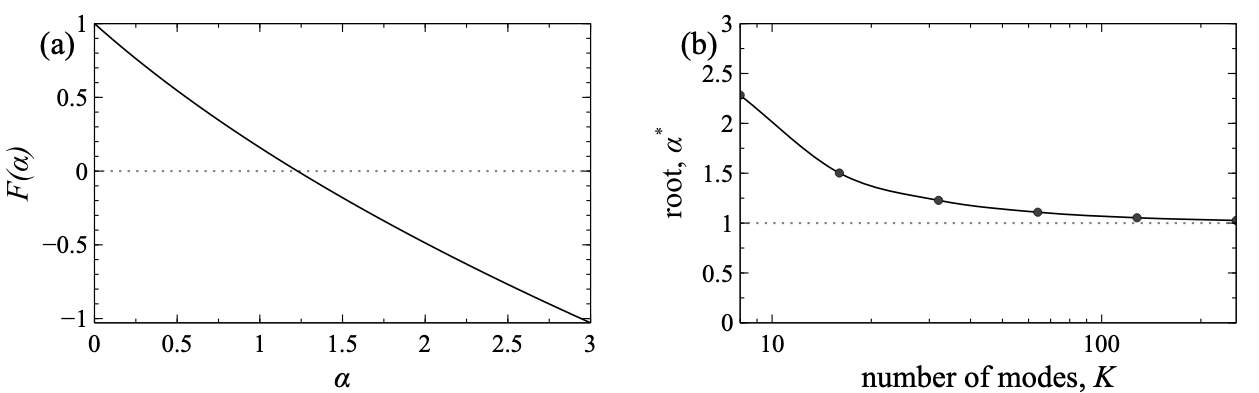}
\caption{
The value $\alpha^*$ used in the approximating measure \cref{dgSig2} is selected as the root of \cref{F_alpha}.
(a) $F(\alpha)$ from \cref{F_alpha} with $\nmodes = 32$ and $\beta' = 20$.
(b) Root $\alpha^*$ for increasing values of $\nmodes$ with $\bprime = 20$.}
\label{alpha_plots}
\end{figure}

\Cref{alpha_plots} provides a visualization of how $\alpha^*$ is selected and how its value depends on the cutoff wavenumber $\nmodes$. \Cref{alpha_plots}(a) shows the function $F(\alpha)$ defined in \cref{F_alpha} for values $\nmodes = 32$ and $\bprime = 20$. The root occurs at approximately $\alpha^* = 1.2$. \Cref{alpha_plots}(b) shows the value of $\alpha^*$, as computed with a numerical root-finder, for the values $\nmodes = 8, 16, 32, 64,$ and 128. Notice that $\alpha^* \to 1$ as $\nmodes \to \infty$, which can be confirmed by directly taking the limit of \cref{F_alpha}.

\subsection{Evaluating the proposal distribution}

Recall that the Gibbs measure \cref{Gibbs1} has an associated probability density, $f:\sphere \to \RR$, supported on the hypersphere. This $f$ is the target distribution we seek to sample, and \cref{f_compute} shows how to compute $f(\xhv)$ for any $\xhv \in \sphere$. Likewise, the approximating measure \cref{dgSig2} has an associated probability density $g:\sphere \to \RR$ supported on the hypersphere and defined by
\begin{equation}
\label{g_defn}
\int_{\sphere} \phi(\xhv) \, g(\xhv) \, d\mu_0(\xhv) = \int_{\RR^{2 \nmodes}} \phi(\xhv) \, d\gamma_{\Sig} (\xv) \, ,
\end{equation}
for any measurable $\phi:\sphere \to \RR$. The density $g$ will be the proposal distribution of the rejection algorithm. Since the underlying measure $d\gamma_{\Sig}$ from \cref{dgSig1} is Gaussian, it is straightforward to sample $g$. One first draws a sample $\Xv$ from $d\gamma_{\Sig}$ and then projects onto $\sphere$ via $\Xhv = \Xv/\abs{\Xv}$. However, to use $g$ as a proposal distribution, it is necessary to compute the ratio $f(\xhv)/g(\xhv)$ to determine the acceptance probability of a given $\xhv \in \sphere$. While \cref{f_compute} shows how to compute $f(\xhv)$, we have not yet shown how to compute $g(\xhv)$.

In a rejection algorithm, it suffices to compute the ratio $f(\xhv)/g(\xhv)$ up to a multiplicative constant. Therefore the value $Z$ in \cref{f_compute} and other constants that arise will not be needed. For this reason, we will use $\Const$ to represent various constants that arise in this section. Even though these constants are represented by the same symbol, they do not necessarily have the same value. 

As seen in \cref{g_defn}, evaluating $g(\xhv)$ for a particular value of $\xhv \in \sphere$ requires integrating the measure $d\gamma_{\Sig}(\xv)$ over all values of $\xv \in \RR^{2\nmodes}$ for which $\xv/\abs{\xv} = \xhv$. These values comprise a ray emanating from the origin and passing through the point $\xhv$, given explicitly by $\xv = r \, \xhv$ for $r\in(0,\infty)$. Evaluating \cref{dgSig1} along this ray gives
\begin{align}
d \gamma_{\Sig} \Big|_{\xv = r \, \xhv} \propto 
\exp \left( -a(\xhv) \, r^2 \right)  \, ,
\end{align}
where the parameter $a(\xhv)$ depends on the direction vector $\xhv \in \sphere$ as given by
\begin{align}
\label{a_defn}
a(\xhv) = \sum_{k=1}^{\nmodes} \frac{\xhat_k^2 + \xhat_{\nmodes+k}^2}{ 2\sigma_k^2 } \, .
\end{align}
Inserting into \cref{g_defn} gives
\begin{equation}
\label{g_integral}
g(\xhv) = \Const \cdot \int_0^{\infty} \exp(-a(\xhv) \, r^2) \, r^{2\nmodes-1} \, dr \, .
\end{equation}
To compute the integral above, we first perform a change of variables with $u = r^2$,
\begin{equation}
\label{u_subs}
\int_0^\infty \exp(-a r^2) \, r^{2\nmodes-1} \, dr 
= \frac{1}{2} \int_0^{\infty} \exp(-a u) \, u^{\nmodes-1} \, du \, .
\end{equation}
We then use a known Laplace transform
\begin{equation}
\label{Laplace}
\mathcal{L}[t^{\nmodes-1}] 
= \int_0^{\infty} \exp(-s t) \, t^{\nmodes-1} \, du = {(\nmodes-1)!} \, {s^{-\nmodes}} \, .
\end{equation}
Using \cref{g_integral,u_subs,Laplace}, the proposal distribution $g(\xhv)$ is then given by
\begin{equation}
g(\xhv) = \Const \cdot a(\xhv)^{-\nmodes} \, .
\end{equation}
Inserting \cref{a_defn} and the standard deviations given by \cref{sigma} then gives
\begin{equation}
\label{g_eqn}
g(\xhv) = \Const \cdot \left( 1 + \frac{\alpha^* \bprime}{\nmodes^3} 
\sum_{k=1}^{\nmodes} k^2 \left( \xhat_k^2 + \xhat_{\nmodes+k}^2 \right) \right)^{-\nmodes} \, .
\end{equation}
This is the formula we will use to evaluate the proposal distribution $g$ on the unit hypersphere. 

We remark that a nice check on the above formula is possible. Taking the limit as $\nmodes \to \infty$ of \cref{g_eqn} and using the limit-definition of the exponential gives
\begin{equation}
\label{g_limit}
g(\xhv) \to \Const \cdot \exp \left( -\frac{\alpha^* \bprime}{\nmodes^2} 
\sum_{k=1}^{\nmodes} k^2 \left( \xhat_k^2 + \xhat_{\nmodes+k}^2 \right) \right)
\qquad \text{as } \nmodes \to \infty \, .
\end{equation}
Inspection of \cref{F_alpha} shows that $\alpha^* \to 1$ in this limit, in which case \cref{g_limit} matches \cref{linearGibbs1}.

\section{The rejection algorithm}
\label{Sec:Algorithms}

We now describe the rejection algorithm that will enable efficient sampling of the Gibbs measure,  \cref{Gibbs1}. Associated with this Gibbs measure is the probability density $f:\sphere \to \RR$ supported on the hypersphere. The density $f$ is straightforward to evaluate through \cref{f_compute}, however, it is not straightforward to sample $f$ directly. We therefore use density $g:\sphere \to \RR$, defined by \cref{g_defn}, as a proposal distribution. Since the measure $d\gamma_{\Sig}$ defining $g$ is Gaussian, it is straightforward to sample $g$; one first draws a sample $\Xv$ from $d\gamma_{\Sig}$ and then projects onto $\sphere$ via $\Xhv = \Xv/\abs{\Xv}$. Further, \cref{conv_result} shows that $g$ converges to $f$ in a particular limit. It is therefore reasonable to expect that $g$ lies fairly close to $f$ in general, at least more so than the naive, uniform distribution.
For these reasons, we expect $g$ to satisfy the two criteria for a good proposal distribution: it approximates $f$ and is easy to sample.
The numerical tests in \cref{Sec:Results} confirm that $g$ yields acceptance rates much higher than the uniform distribution, even well outside the strict validity of \cref{conv_result}.

We will employ a standard rejection algorithm; that is, we first draw a sample from the proposal distribution $g(\xhv)$ and then accept the sample with probability proportional to $f(\xhv)/g(\xhv)$. Repeating this process yields a collection of samples from the target distribution $f(\xhv)$ \cite{Ross2010}. Alternatively, one could employ the sampling-importance resampling (SIR) algorithm with $g$ as the proposal distribution \cite{skare2003improved}, although we do not implement it here.

\subsection{Computing the acceptance probability and rejection constant}

In the proposed rejection algorithm, the acceptance probability is proportional to the ratio $f/g$. For simplicity, we take $Const$ in \cref{g_eqn} to equal $Z^{-1}$ so that constants cancel to give
\begin{equation}
\label{fg_ratio}
\frac{f(\xhv)}{g(\xhv)} =
\exp\left( -\beta \, \Ham_{\nmodes}[ \Ufun( \xhv ) ] \right)
\left( 1 + \frac{\alpha^* \bprime}{\nmodes^3} 
\sum_{k=1}^{\nmodes} k^2 \left( \xhat_k^2 + \xhat_{\nmodes+k}^2 \right) \right)^{\nmodes} \, .
\end{equation}

In order to compute $\Ham_{\nmodes}$ above it is necessary to compute $\Ham_2$ and $\Ham_3$. $\Ham_2$ can be efficiently computed with \cref{H2Pars}. For $\Ham_3$, Abramov et al.~\cite{Abramov2003} derived the triple-summation formula
\begin{align}
\Ham_3 &= \frac{\pi}{3} \sum_{\substack{ k_1 + k_2 + k_3 = 0 \\ \abs{k_1}, \abs{k_2}, \abs{k_3} \le \nmodes}}\uhat_{k_1} \uhat_{k_2} \uhat_{k_3} \, .
\end{align}
This formula can be rearranged as a double summation
\begin{align}
\label{H3_double}
\Ham_3 &= 2\pi \sum_{n=1}^{\nmodes}
\, Re \left( \uhat_{n}^* \sum_{k=1}^{n-1} \uhat_k \uhat_{n-k} \right) \, .
\end{align}
Therefore, $\Ham_3$ can be computed with $O(\nmodes^2)$ operations. We remark that, alternatively, $\Ham_3$ could be computed from the definition \cref{H3H2} using the FFT with $O(\nmodes \log \nmodes)$ operations. However, for moderate values of $\nmodes$ in the range $8 \le \nmodes \le 64$, \cref{H3_double} is typically more efficient and hence the form that we will use to compute $\Ham_3$.

In practice, we reduce the number of control parameters in \cref{fg_ratio} by recasting in terms of $\bprime$ and the ratio $C_3/C_2$. That is, we use \cref{TKdVHam} and \cref{beta_prime} to rewrite
\begin{align}
\label{beta_Ham}
\beta \, \Ham_{\nmodes} = 
\frac{\bprime}{\En_0 \nmodes^2} \left( \Ham_2 - \frac{C_3}{C_2} \Ham_3 \right) \, .
\end{align}
Therefore, the only parameters that must be specified are $\nmodes, \En_0, \bprime$, and $C_3/C_2$; It is not necessary to specify $\beta$, $C_2$, and $C_3$ individually.
In \cref{appA}, we provide exact solutions for $\Ham_2$ and $\Ham_3$ in the case of $\nmodes=2$. These solutions are used to validate our numerical computation of $\Ham_2$ and $\Ham_3$.

We now have everything needed to compute the ratio $f/g$ through \cref{fg_ratio}. However, it remains to compute the rejection constant. That is, we seek a constant $M$ to ensure that
\begin{equation}
\frac{f(\xhv)}{M g(\xhv)} \le 1
\quad \text{for all } \xhv \in \sphere \, ,
\end{equation}
so that the acceptance probability is no greater than one. The optimal such constant is given by
\begin{equation}
\label{M_sup}
M = \sup_{\xhv \in \sphere} \frac{f(\xhv)}{g(\xhv)}  \, .
\end{equation}
Since the unit hypersphere $\sphere$ is a compact set, the supremum is equal to the maximum value. We can therefore solve a numerical optimization problem to determine $M$. Since the dimension of $\sphere$ is not necessarily small (e.g.~the typical range $8 \le \nmodes \le 32$ yields $15 \le \dim(\sphere) \le 63$), it is important to have a good initial guess for the numerical optimization.

To obtain a good initial guess, we consider the derivation of the proposal distribution $g$ and the physical meaning of the ratio $f/g$. In particular, $g$ was chosen specifically to account for the presence of $\Ham_2$ in the Hamiltonian, whereas $f$ includes contributions from both $\Ham_2$ and $\Ham_3$. In particular, since $f$ is proportional to $\exp(-\beta C_2 \Ham_2 + \beta C_3 \Ham_3)$ and $g$ is approximately proportional to $\exp(-\beta C_2 \Ham_2)$, the $\Ham_2$ contributions approximately cancel in the ratio, giving
\begin{equation}
\frac{f(\xhv)}{g(\xhv)} \approx
\Const \cdot \exp(\beta C_3 \Ham_3) \, .
\end{equation}
Thus maximizing $\Ham_3$ should provide a good starting guess for maximizing the ratio $f/g$.

\begin{figure}
\centering
\includegraphics[width = \linewidth]{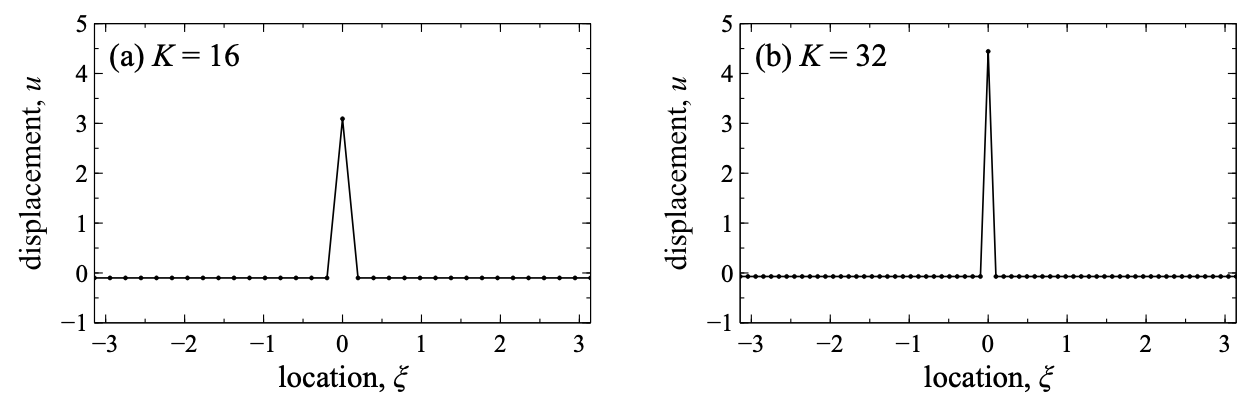}
\caption{
Visualization of the Dirichlet kernel for $\nmodes = $ 16 and 32. This function is used as the initial guess to numerically optimize the ratio $f/g$ on the unit hypersphere.
}
\label{Dirichlet}
\end{figure}

Recall from \cref{H3H2}, that $\Ham_3$ measures the skewness, $\int_{-\pi}^\pi u^3 d\xc$, of the sampled surface displacement. We seek to maximize this quantity within the space of functions having Fourier modes up to wavenumber $\nmodes$ and satisfying the constraints \cref{MomEn} of zero momentum and fixed energy. In physical space, the maximum is achieved by a Dirac-delta approximation, that is, a spike-like function with a single, large positive value. This function must be shifted down slightly to satisfy the zero-mean condition, $\int_{-\pi}^\pi u d\xc = 0$. The best finite-mode Dirac-delta approximation is known as the Dirichlet kernel \cite{Strauss}, characterized by a constant spectrum. 
In spectral space, the zero-mean Dirichlet kernel with peak centered at $\xi = 0$ is given by
\begin{equation}
\label{optim_guess}
\uhat_k = \sqrt{ \frac{\En_0}{2 \pi \nmodes} }
\qquad \text{for } k = \pm 1,2,\cdots, \nmodes \, .
\end{equation}
\Cref{Dirichlet} illustrates this Dirichlet kernel in the cases $\nmodes = $ 16 and 32. Both show a large positive value at $\xi = 0$ and a small negative value for $\xi \ne 0$, where these two values are chosen to satisfy zero mean and unit energy.

Using \cref{optim_guess} as the initial guess, we numerically optimize the ratio $f/g$ with the Nelder-Mead algorithm \cite{lagarias1998convergence} in order to determine the rejection constant $M$.
We find that the numerical optimization typically converges in only a few iterations, confirming that \cref{optim_guess} provides an effective starting guess.

\subsection{The rejection algorithm}

All of the pieces are now in place to implement the rejection algorithm. To summarize, the algorithm is as follows:
\begin{enumerate}
\item Precompute the the rejection constant $M$ by numerically maximizing \cref{M_sup} using \cref{optim_guess} as the initial guess.
\item Draw a sample $\Xv \in \RR^{2\nmodes}$ from the anisotropic Gaussian distribution $d\gamma_{\Sig}$ given by \cref{sigma,dgSig1} with $\alpha = \alpha^*$ selected as the root of \cref{F_alpha}.
\item Project the sample $\Xv$ onto $\sphere$ via $\Xhv = \Xv/\abs{\Xv}$. By doing so, $\Xhv$ represents a random sample drawn from density $g$.
\item For the sample $\Xhv$, compute the acceptance probability $f(\Xhv)/(M g(\Xhv))$ using \cref{fg_ratio}, with \cref{H2Pars,H3_double} used to compute $\Ham_2$ and $\Ham_3$. Decide to accept the sample or not based on the computed probability. If accepted, $\Xhv$ represents a random sample drawn from density $f$.
\edit{Through \cref{uhxh}, an accepted $\Xhv$ produces a Fourier-coefficient set, $\uhat_k$ for $k=1,\cdots,\nmodes$, of an accepted wave field.}
\item Go to step 2 and repeat until a desired number of samples have been accepted.
\end{enumerate}

We note that it is trivial to parallelize the above algorithm. In particular, step 1 should be precomputed and then steps 2--5 can be performed in parallel with no communication required between computational nodes. In the following section, we report numerical results on both serial and parallelized versions of this algorithm.


\section{Results} 
\label{Sec:Results}

We now present numerical results of the novel rejection algorithm. To provide physical context, we first visualize some wave fields and statistical data generated by the algorithm in different parameter regimes of TKdV. We then evaluate the performance of the algorithm, including: (i) a comparison against a naive version that uses the \edit{spectrally} uniform measure as a proposal distribution, and (ii) results on the parallel speedup of the algorithm. We then examine characteristics of individual wave fields sampled by the algorithm, including extreme events.

\subsection{Physical wave fields generated by the  algorithm}

The main purpose of the new rejection algorithm is to efficiently generate independent, random samples of wave fields arising from the TKdV system without directly simulating the dynamics. To orient the reader to the physics of interest, we begin by showing a few example wave fields generated by the algorithm in different parameter regimes.
In all of these examples, we fix the maximum wavenumber $\nmodes = 16$ and the energy $\En_0 = 1$, while varying the normalized inverse temperature $\bprime$ and the nonlinearity-to-dispersion ratio $C_3/C_2$. Recall that due to \cref{beta_Ham}, these are all of the control parameters to be specified (i.e.~it is not necessary to specify $\beta, C_2$, and $C_3$ individually).
\edit{One feature we analyze below is the distribution of surface displacements that results from the Gibbs measure. We remind the reader that the Gibbs measure exists in spectral space, i.e.~it is a distribution of the wave field's Fourier coefficients, whereas the  surface-displacement distribution exists in physical space. Gaussianity (or non-Gaussianity) of surface displacement statistics does not imply Gaussianity (or non-Gaussianity) of the underlying Gibbs measure. For example, the spectrally uniform Gibbs measure $d\mu_0$, that arises in the case $\bprime=0$, produces Gaussian statistics of the surface displacement and not a uniform distribution of surface displacements \cite{SunMoore2023}.}

\begin{figure}
\centering
\includegraphics[width=\linewidth]{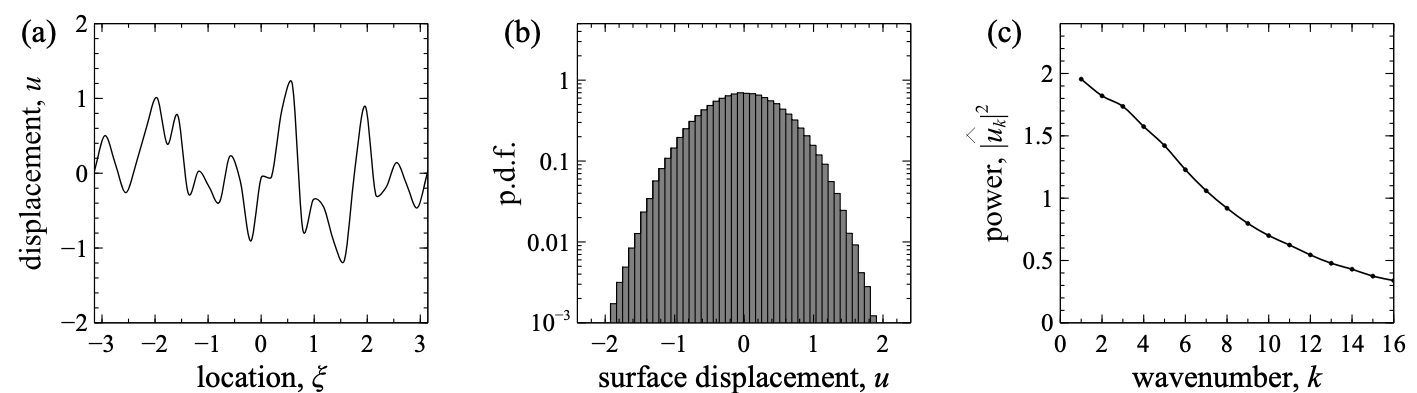}
\caption{
Visualizing wavefield statistics in the case of {\em} linear TKdV ($C_3/C_2 = 0$) with $\bprime = 40$. (a) A representative wave field sampled by the rejection algorithm. (b) Histogram of surface displacements over 5,000 samples. The distribution is symmetric and nearly Gaussian. (c) Power spectrum averaged over all 5,000 samples. The spectrum decays gradually with wavenumber.
}
\label{wave1}
\end{figure}

To begin, \cref{wave1} displays wave information for the case of {\em linear} TKdV, i.e.~$C_3/C_2 = 0$, with a moderate inverse temperature $\bprime = 40$. We use the rejection algorithm to sample the Gibbs measure until 5,000 samples are accepted.
\Cref{wave1}(a) shows a representative wave field, \edit{obtained by performing a Fourier transform of a spectrum sampled from the Gibbs measure}. The wave exhibits positive and negative surface displacements of roughly equal magnitudes. This plot represents an ordinary wave field with no extreme values.
\Cref{wave1}(b) combines all 5,000 samples into a histogram of surface displacements. Given the logarithmic vertical axis, the roughly parabolic shape implies a nearly normal distribution. That is, this set of parameters produces {\em Gaussian statistics of the surface displacement}, as is typically observed for {linear} TKdV \cite{Majda2019, MooreJNS2020, SunMoore2023}. \edit{Note that the Gibbs distribution itself, which exists in spectral space, is not Gaussian.}
\Cref{wave1}(c) shows the power spectrum averaged over all 5,000 samples. The gradual decay of the spectrum is due to $\Ham_2$, which penalizes higher wavenumbers as seen through \cref{H2Pars}.

\begin{figure}
\centering
\includegraphics[width=\linewidth]{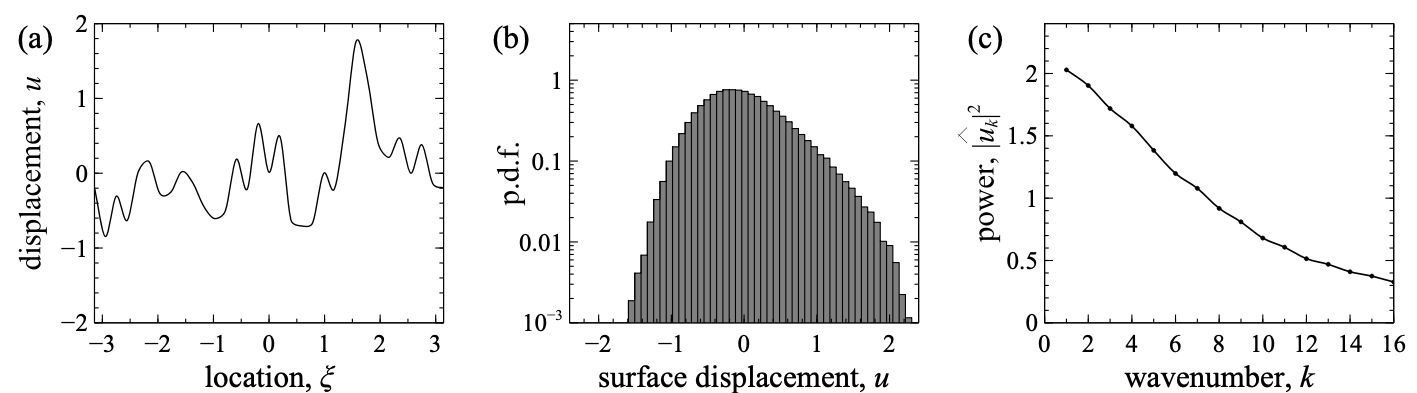}
\caption{
The effects of strong nonlinearity, $C_3/C_2 = 120$, with $\bprime = 40$ as before. (a) The representative wave field features a prominent peak occurring near $\xi = 1.5$. (b) The histogram exhibits  positive skewness, indicating that large, positive displacements are consistently favored in the ensemble of 5,000 accepted samples. (c) The spectrum decays gradually, similar  to the previous case.
}
\label{wave2}
\end{figure}

\Cref{wave2} shows the case of stronger nonlinearity, $C_3/C_2 = 120$, with the other parameters held fixed. \Cref{wave2}(a) shows the representative wave field, distinguished in this case by a prominent peak near $\xi = 1.5$. This wave field is not simply an outlier. Indeed, the histogram of 5,000 accepted samples in \cref{wave2}(b) exhibits a strong positive skewness (skewness = 0.59), indicating that large, positive values of $u$ are more likely in the ensemble.
The reason for this behavior is that the stronger nonlinearity (i.e.~larger $C_3/C_2$) increases the relative contribution of $\Ham_3 \propto \int_{-\pi}^\pi u^3 d\xc$, which favors positive skewness. \Cref{wave2}(c) shows that the spectrum decays gradually, as before. To leading order, the spectrum is unaffected by nonlinearity \cite{SunMoore2023}.

These two cases illustrate some general principles underlying TKdV wave statistics. In particular, for $\bprime > 0$, the $\Ham_2$ component penalizes high frequencies in the Gibbs ensemble, creating a gradually decaying power spectrum. If nonlinearity is absent, the surface-displacement distribution is symmetric \cite{SunMoore2023}, as seen in \cref{wave1}(b). If nonlinearity is present, the $\Ham_3$ component favors positively skewed wave fields, leading to skewness in the ensemble histogram (\cref{wave2}(b)) and increasing the likelihood of large, positive displacements. In particular, nonlinearity increases the likelihood of extreme events, as illustrated by the large peak seen in \cref{wave2}(a). In the more trivial case of $\bprime = 0$ (not shown), the Gibbs measure is independent of the Hamiltonian, which results in nearly Gaussian surface-displacement statistics \cite{SunMoore2023}.

\subsection{Performance of the algorithm}
\label{Performance}

We now evaluate the performance of the rejection algorithm equipped with the new proposal distribution $g$, defined in \cref{g_defn}. In particular, the average acceptance rate determines the speed at which random wave fields can be generated. Therefore, comparing the average acceptance rate achieved by $g$ versus the uniform measure, \cref{mu0}, indicates the degree of performance improvement. To simplify nomenclature, we refer to the use of $g$ as the `improved rejection algorithm' and the use of the uniform measure as the `naive' version.

\begin{table}
\begin{center}
\begin{tabular}{ |c|c|c|c|c| }
\hline \multicolumn{5}{|c|}{Performance Test, $\nmodes = 16$} \\
\hline
$\bprime$ & $C_3/C_2$ & Skewness & Acceptance Rate & Improvement \\ \hline
   &0   &0     &$0.95$              &430 \\
   &60  &0.11  &$2.4\times 10^{-2}$ &60 \\
20 &120 &0.23  &$6.9\times 10^{-4}$ &25 \\ 
   &180 &0.37  &$2.1\times 10^{-5}$ &16 \\ 
   &240 &0.53  &$7.7\times 10^{-7}$ &12 \\ 
   &300 &0.73  &$3.2\times 10^{-8}$ &11 \\ 
\hline
   &0   &0     &0.88                &49,000 \\
   &30  &0.13  &$2.3\times 10^{-2}$ &3,400 \\
40 &60  &0.25  &$6.7\times 10^{-4}$ &690 \\ 
   &90  &0.40  &$2.2\times 10^{-5}$ &280 \\ 
   &120 &0.59 &$8.0\times 10^{-7}$ &130 \\ 
   &150 &0.81  &$3.6\times 10^{-8}$ &74 \\ 
\hline
   &0	 &0    &0.80                &1,600,000 \\
   &20  &0.14 &$2.3\times 10^{-2}$ &120,000 \\ 
60 &40  &0.27 &$7.0\times 10^{-4}$ &23,000 \\
   &60  &0.42 &$2.4\times 10^{-5}$ &4,300 \\
   &80  &0.66 &$9.2\times 10^{-7}$ &1,200 \\
   &100 &0.85 &$4.6\times 10^{-8}$ &590 \\
\hline
\end{tabular}
\caption{Performance of the improved rejection algorithm for $\nmodes = 16$. Columns 1--2: the parameter values $\bprime$ and $C_3/C_2$ chosen for the test; $\En_0 = 1$ in all cases. Column 3: the skewness of the resulting ensemble, which serves as a proxy for anomalous wave behavior. Column 4: the average acceptance rate of the improved algorithm. The acceptance rate is extremely high for linear TKdV, and it decreases as nonlinearity, $C_3/C_2$, increases. Column 5: the factor by which the novel proposal distribution improves the acceptance rate over the uniform distribution. The test shows an improvement factor of 1--6 orders of magnitude.
}
\label{Table1}
\end{center}
\end{table}

\Cref{Table1} summarizes results from a series of tests with cutoff wavenumber $\nmodes = 16$, energy $\En_0 = 1$, and three inverse temperatures, $\bprime = $ 20, 40, and 60. For each $\bprime$, we increment the nonlinearity parameter, $C_3/C_2$, to obtain the same increments of $\bprime C_3/C_2$, since this combination is the prefactor of $\Ham_3$ in the Hamiltonian. As seen in \cref{H3H2}, $\Ham_3$ controls the skewness of the wave field, and so equal increments of $\bprime C_3/C_2$ results in nearly equal increments of ensemble skewness. In this way, the numerical test probes a similar physical regime for each $\bprime$.
The second column of \Cref{Table1} shows the ensemble skewness measured in the numerical tests and confirms that, for each $\bprime$, the skewness increments by roughly the same amount as nonlinearity increases.

Column 3 of \cref{Table1} shows acceptance rates achieved by the improved rejection algorithm, while Column 4 shows how many times greater this acceptance rate is than that achieved by the naive version. In all cases, the improvement factor is much greater than one, and in most cases by several orders of magnitude. These results indicate that $g$ vastly outperforms the uniform measure as a proposal distribution, \edit{increasing the number of accepted samples by 1--6 orders of magnitude.}

More specifically, the improved rejection algorithm produces extremely high acceptance rates when nonlinearity is small. For $C_3/C_2 = 0$, the acceptance rates are all above 80\%, indicating that the proposal distribution $g$ approximates the true Gibbs density, \cref{f_compute}, extremely well. In these cases, the performance gap between the improved and naive proposal distributions is the largest. In the most extreme case of $\bprime = 60$, {the improved algorithm produces 1.6 million times more accepted samples than does the naive version.}

As nonlinearity, $C_3/C_2$, increases, the acceptance rate of the improved algorithm decreases, but it still outperforms the naive version by a wide margin. The most physically interesting cases are those of strong skewness (roughly skewness $> 0.4$) indicative of highly anomalous statistics and an increased likelihood of extreme events. In this regime, the improved algorithm achieves acceptance rates in the range $10^{-5}$ to $10^{-8}$. Interestingly, the improvement factor over the naive version increases with $\bprime$.
That is, in the regime of high skewness, the novel proposal distribution improves performance by an order of magnitude for $\bprime = 20$ but three orders of magnitude for $\bprime = 60$.

\begin{table}
\begin{center}
\begin{tabular}{ |c|c|c|c|c| }
\hline \multicolumn{5}{|c|}{Performance Test, $\nmodes = 32$} \\
\hline
$\bprime$ & $C_3/C_2$ & Skewness & Acceptance Rate & Improvement \\ \hline
   &0   &0    &0.99                 &360 \\
   &150 &0.04 &$3.4 \times 10^{-2}$ &13 \\
   &300 &0.08 &$1.2 \times 10^{-3}$ &41 \\ 
20 &450 &0.11 &$4.2 \times 10^{-5}$ &27 \\
   &600 &0.15 &$1.5 \times 10^{-6}$ &20 \\
   &750 &0.19 &$5.5 \times 10^{-8}$ &14 \\
   &900 &0.25 &$2.1 \times 10^{-9}$ &14 \\
\hline
   &0   &0  &0.96                    &66,000 \\
   &75  &0.05  &$3.3 \times 10^{-2}$ &1,600 \\ 
   &150 &0.07  &$1.1 \times 10^{-3}$ &79 \\
40 &225 &0.12  &$4.1 \times 10^{-5}$ &2 \\
   &300 &0.16  &$1.4 \times 10^{-6}$ &700 \\
   &375 &0.22  &$5.5 \times 10^{-8}$ &410 \\
   &450 &0.26  &$2.1 \times 10^{-9}$ &$330^*$ \\
\hline
   &0   &0    &0.92                 &2,500,000 \\
   &50  &0.04 &$3.1 \times 10^{-2}$ &71,000 \\
   &100 &0.08 &$1.1 \times 10^{-3}$ &1,900 \\
60 &150 &0.13 &$3.9 \times 10^{-5}$ &390 \\
   &200 &0.17 &$1.4 \times 10^{-6}$ &13 \\
   &250 &0.23 &$5.5 \times 10^{-8}$ &2.3 \\
   &300 &0.26 &$2.0 \times 10^{-9}$ &5,700$^*$ \\
\hline
\end{tabular}
\caption{Performance of the improved rejection algorithm as compared to the naive version for $\nmodes = 32$. As before, the improved proposal distribution increases the number of accepted samples by several orders of magnitude. To measure the improvement ratio we require at least 2,000 samples from the improved algorithm and at least 100 from the naive version. The asterisk indicates cases in which the naive version was unable to produce 100 samples due to an extremely low acceptance rate.
}
\label{Table2}
\end{center}
\end{table}

\Cref{Table2} summarizes results of a similar test with larger cutoff wavenumber $\nmodes = 32$. As before, we increment $C_3/C_2$ so as to obtain the same increments in $\bprime C_3/C_2$ and thus similar increments in skewness. When nonlinearity is absent ($C_3/C_2=0$) the improved algorithm still enjoys extremely high acceptance rates ($>90\%$), indicating that the proposal distribution continues to approximate the Gibbs measure extremely well. The improvement over the uniform measure is extremely large, up to 2.5 million in the case $\bprime = 60$. As nonlinearity increases, the ensemble skewness increases and the acceptance rate decreases. The improved algorithm still enjoys much better performance than the naive version in this regime, increasing the number of accepted samples by 1-5 orders of magnitude. 

A few surprises emerge in \cref{Table2}. First, in the case $\bprime=40$, the improvement factor depends non-monotonically on $C_3/C_2$; it initially decreases with increasing $C_3/C_2$, reaches a minimum at $C_3/C_2 = 225$, and then increases. We have verified this non-monotone behavior is repeatable by rerunning the tests with different random seeds, but we do not have a simple explanation for it. Second, comparing cases of large $C_3/C_2$ across the different values of $\bprime$, it seems that the improved algorithm enjoys a sweet spot at $\bprime =40$, where its advantage over the naive version is maximal. 
For the case $\bprime C_3/C_2 = $15,000 (second to largest value of $C_3/C_2$ in each group), the improvement factor is 410 for $\bprime = 40$, but only 14 and 2.3 for $\bprime =$ 20 and 60 respectively.

We note that to measure the improvement ratio, we required at least 2,000 samples to be drawn by the improved algorithm and at least 100 by the naive algorithm. In the two cases $(\bprime, C_3/C_2) = (40, 450)$ and $(60,300)$ the acceptance rate of the naive version was so low that the test could not be completed over the course of weeks. We therefore estimated the naive acceptance rate and the corresponding improvement ratio from the low number of samples that were drawn. These values are indicated with an asterisk in \cref{Table2}.

\subsection{Parallel speedup of the algorithm}

\begin{figure}
\centering
\includegraphics[width = 0.6 \linewidth]{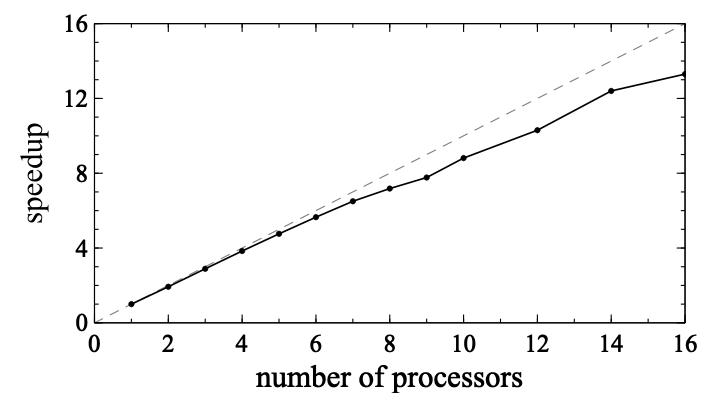}
\caption{
The parallel speedup of the rejection algorithm versus the number of processors. Test conducted on an Apple M2 Ultra chip with 16 performance cores. The cutoff wavenumber is $\nmodes = 128$ and the number of nominal samples per thread is $10^4$.
}
\label{speedup}
\end{figure}

A significant advantage of a rejection algorithm is the favorable parallelization properties. When equipped with a good proposal distribution and run in parallel, a simple rejection algorithm can outperform more sophisticated methods like MCMC. To demonstrate the favorable parallelization properties, we run the improved rejection algorithm in parallel with up to 16 processors. \Cref{speedup} shows the speedup achieved versus the number of processors. As seen in the figure, the algorithm achieves nearly optimal speedup for up to 8 processors. Beyond 8 processors, the speedup diminishes as would be expected for any algorithm due to the required communication between processors and memory. The reduction in speedup, however, is small. For the largest case of 16 processors, the speedup is still 80\% of the optimal value.

\subsection{Sampling individual wave fields}

A major advantage of a rejection algorithm over importance sampling is the ability to draw individual samples from a distribution, rather than merely evaluating statistical properties of the distribution (e.g. computing means or generating histograms). To demonstrate this capability, we now examine the characteristics of individual samples drawn from the Gibbs measure, with particular attention paid to the most extreme waves.

\begin{figure}
\centering
\includegraphics[width = 1.0 \linewidth]{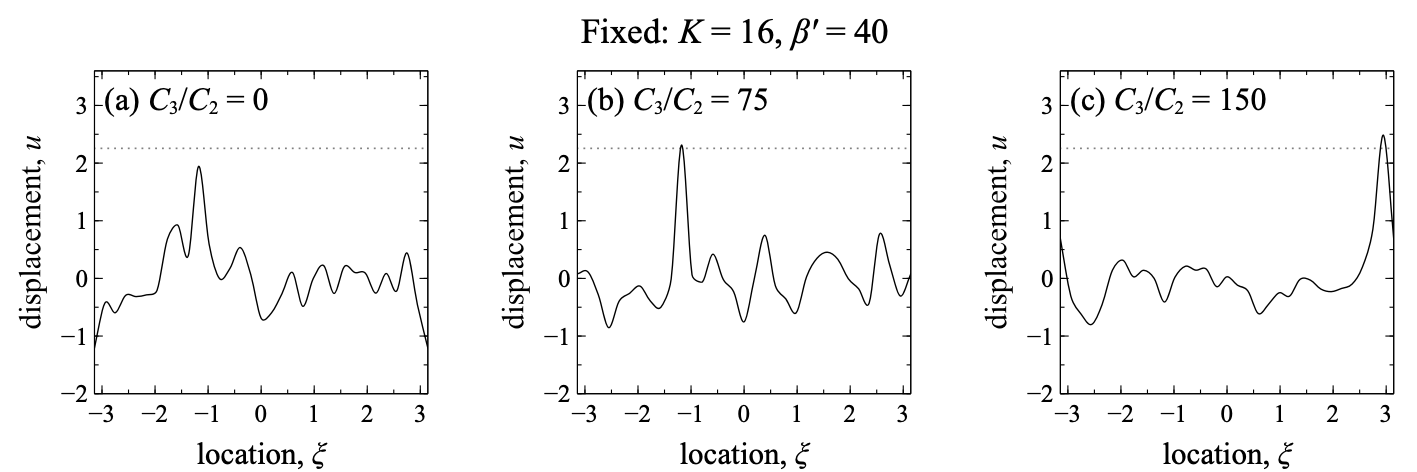}
\caption{
Sampling individual wave fields to search for extreme events. For each $C_3/C_2$, we show the individual wave field drawn from an ensemble of 500 that possesses the largest surface displacement, i.e.~a 1-in-500 event. Increasing nonlinearity (left to right) increases the value of the maximum. In the last two cases, the maximum displacement exceeds the $4\sigma$ threshold (faint dotted line). The parameters $\nmodes = 16$, $\En_0 = 1$, and $\bprime = 40$ are fixed.
}
\label{big1}
\end{figure}

\begin{figure}
\centering
\includegraphics[width = 1.0 \linewidth]{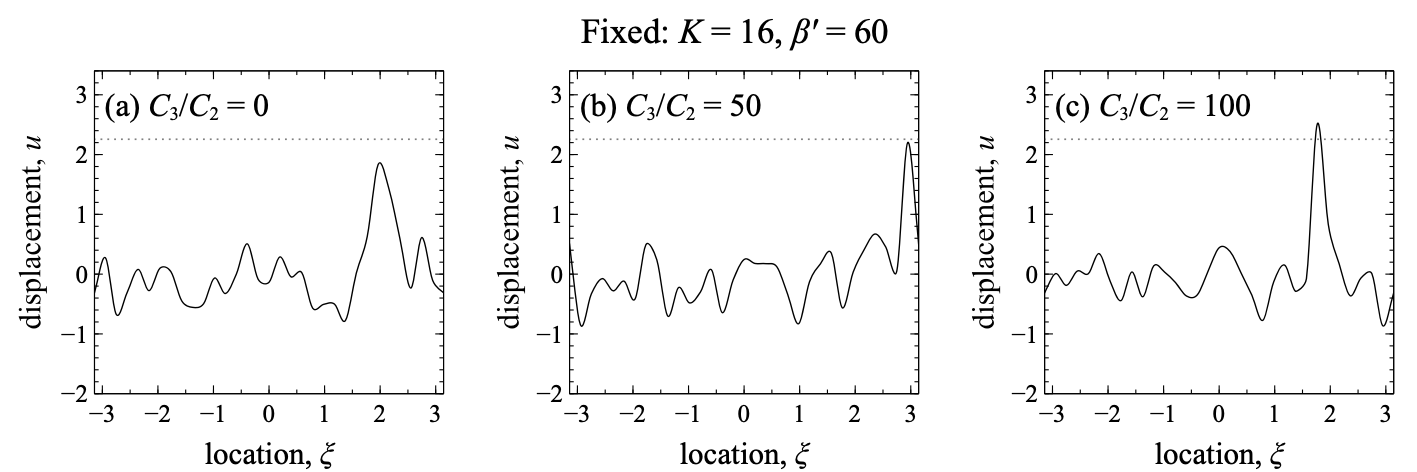}
\caption{
Same test as in \cref{big1} except with $\bprime = 60$ and different increments in $C_3/C_2$. The first two cases of $C_3/C_2$ produce a 1-in-500 event that falls short of the $4\sigma$ threshold. In the last case, the 1-in-500 wave exceeds this threshold.
}
\label{big2}
\end{figure}

\Cref{big1} shows the first test with cutoff wavenumber $\nmodes=16$ and inverse temperature $\bprime = 40$. We consider three values of the nonlinearity ratio, $C_3/C_2 =$ 0, 75, 150. For each, we use the improved rejection algorithm to draw 500 individual wave fields from the Gibbs measure. We then extract the wave field with the largest positive displacement value to display in \cref{big1}; that is, each wave field shown in \cref{big1} represents a 1-in-500 event. For reference, the faint dotted line shows a displacement equal to four standard deviations from the mean, $4\sigma = 2.26$ since $\sigma = \sqrt{E_0/\pi}$ and $E_0=1$. For a Gaussian distribution, roughly one in 16,000 samples lies at least four standard deviations from the mean, and this is sometimes taken to define a rogue-wave event \cite{Bolles2019}. 
We note that each wave field contains $2\nmodes = 32$ surface-displacement values and so a crude estimate for the number of wavefield samples required to achieve the $4\sigma$ threshold under Gaussian statistics is $16,000/32 = 500$.
\Cref{big1}(a) shows that in the linear case, $C_3/C_2 = 0$, the largest surface displacement lies below the $4\sigma$ threshold. \Cref{big1}(b)--(c) show that as nonlinearity increases, the 1-in-500 wave exceeds the $4\sigma$ threshold, with the largest value, $C_3/C_2 = 150$, giving a surface displacement $11\%$ above the threshold. 
To appreciate this example from a different light, we note that, unlike a Gaussian distribution, the maximum displacement of a finite-mode wave field with fixed energy is bounded above. In particular, the Dirichlet kernel achieves the maximum possible displacement (see \Cref{Dirichlet}). For $\nmodes=16$ and $\En_0=1$, the maximum possible displacement is $u_{\text{max}} = 3.09$, and the displacement seen in \cref{big1}(c) is $81\%$ of this value.

\Cref{big2} shows a similar test with the same cutoff wavenumber $\nmodes = 16$, but a larger inverse temperature $\bprime = 60$. The values of the nonlinearity parameter, $C_3/C_2 = 0, 50, 100$, are chosen to yield the same increments in $\bprime C_3/C_2$. As seen in \Cref{big2}, the first two values of $C_3/C_2$ yields a 1-in-500 event that does not meet the $4\sigma$ threshold. \cref{big2}(c) shows that the last case does exceed the threshold by the same $11\%$ observed in \cref{big1}(c). Once again, the fixed energy and cutoff wavenumber imposes a hard cap of $u_{\text{max}} = 3.09$ on the displacement. Our tests suggest that it is difficult to surpass roughly 80\% of this cap, even with significant skewness present in the ensemble (skewness = 0.85 for the ensemble of \cref{big2}(c)).

\begin{figure}
\centering
\includegraphics[width = 1.0 \linewidth]{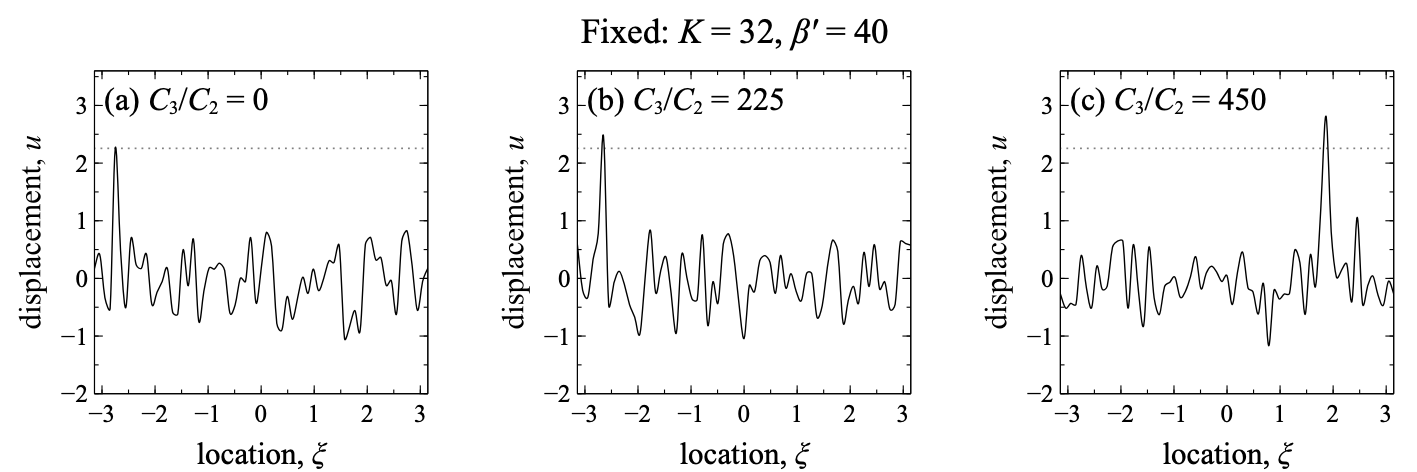}
\caption{
Sampling of individual wave fields with $\nmodes = 32$. For each case, we show the wave field out of 500 that achieves the largest positive displacement. (b)--(c) The latter two cases with strong nonlinearity produce a displacement that exceeds the $4\sigma$ threshold. Fixed parameters are $\bprime=40$ and $\En_0 = 1$.
}
\label{big3}
\end{figure}

\begin{figure}
\centering
\includegraphics[width = 1.0 \linewidth]{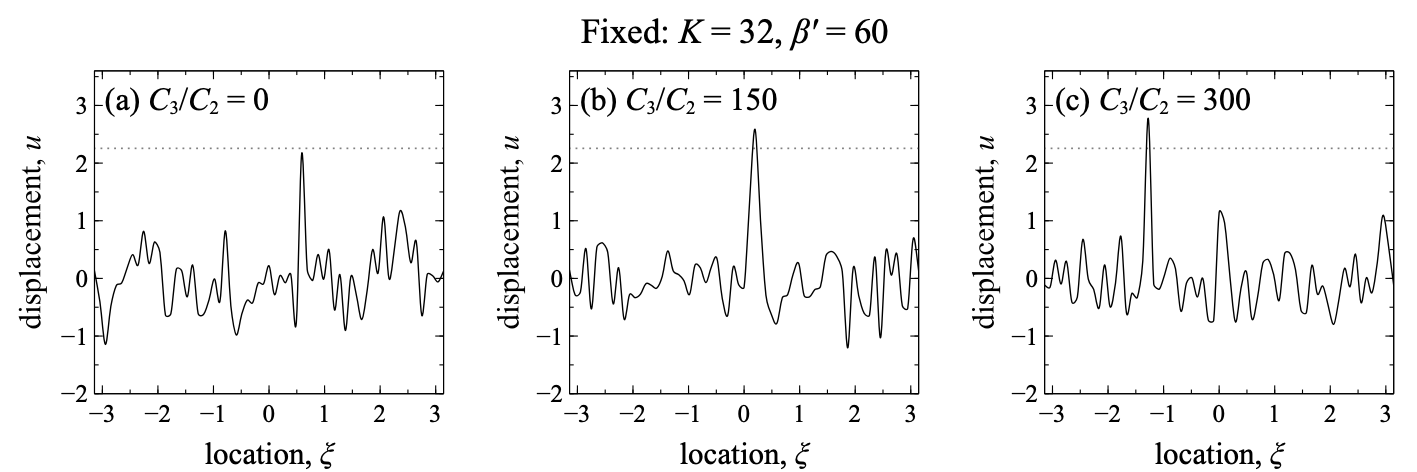}
\caption{
Same test as in \cref{big3} except with larger inverse temperature $\bprime = 60$. (c) The most extreme case $C_3/C_2=300$ produces a displacement that exceeds $4\sigma$ by 24\%, which is the widest margin observed in all tests.
}
\label{big4}
\end{figure}

By increasing the cutoff wavenumber,  we can increase the hard cap on the maximum surface displacement. In \cref{big3,big4} we examine wave fields sampled with $\nmodes = 32$, for which the Dirichlet kernel achieves a maximum possible displacement of $u_{\text{max}} = 4.44$. \Cref{big3} shows the 1-in-500 events for $\bprime = 40$ and three values of $C_3/C_2$. As before, the linear case $C_3/C_2 = 0$ produces a maximal displacement slightly below the $4\sigma$ threshold, while the two nonlinear cases surpass the threshold. In the case of highest nonlinearity, the 1-in-500 maximal displacement exceeds the $4\sigma$ threshold by 24\%. This value is significantly greater than the 11\% achieved in the tests with $\nmodes=16$, despite significantly lower ensemble skewness (0.25 vs.~0.85). This test demonstrates the significance of increasing the hard cap, $u_{\text{max}}$, by increasing the cutoff wavenumber. Unfortunately, achieving comparable values of ensemble skewness with higher $\nmodes$ requires greater computational resources, as seen in \cref{Table1,Table2}.

\Cref{big4} shows a similar test with $\nmodes = 32$ but with a larger inverse temperature of $\bprime = 60$ and the increments of $C_3/C_2$ chosen accordingly. Once again, the linear case produces a 1-in-500 event that falls short of the $4\sigma$ threshold, while both nonlinear cases exceed the threshold. The most nonlinear case, \cref{big4}(c), exceeds the threshold by roughly the same 24\% seen in \cref{big3}(c).

\section{Discussion}
\label{Sec:Discussion}

This paper introduces a novel proposal distribution for generating independent, random samples from the TKdV Gibbs measure. The proposal distribution stems from the analysis of Sun \& Moore (2023) \cite{SunMoore2023}, who showed that it converges to the Gibbs measure in the infinite-mode limit of linear TKdV. Roughly speaking, this proposal distribution  captures the spectral decay of sampled wave fields accurately so that there remains only one criterion for acceptance: sufficient skewness. Numerical tests indicate that the novel proposal distribution produces much higher acceptance rates than the \edit{spectrally} uniform distribution, typically by 1--6 orders of magnitude, even for parameters well outside of the regime in which convergence can be established rigorously. Advantages of the algorithm over more sophisticated methods, like MCMC, include favorable parallelization properties and generation of uncorrelated samples.

Numerical tests demonstrate the algorithm's capability to generate extreme displacements over broad parameter ranges, with several sampled wave fields exceeding the $4\sigma$ criterion of a rogue wave. With a cutoff wavenumber of $\nmodes=16$, tests indicate that $4\sigma$ events can only be generated when the ensemble skewness is sufficiently high (skewness $>0.5$). With $\nmodes =32$, $4\sigma$ events can be generated with lower ensemble skewness (skewness $>0.1$), although the required computational resources are similar.

The algorithmic foundations established here open many exciting avenues for future research. First, the sampling algorithm will be used to build a database of wavefield ensembles corresponding to different regions of parameter space. The database will enable one to efficiently search across parameter space for extreme-wave events and characterize in detail the extreme events that occur in different parameter regimes. The database will also provide a repository of wave fields to be used as initial conditions in dynamic simulations. In particular, running extreme-wave initial conditions {\em backwards} in time will reveal the wave characteristics that precede rogue events. When combined with neural-network feature detection \cite{zhang2018featurenet, pelt2018mixed, qi2024unambiguous}, this process could identify warning signs of impending extreme wave activity.
Lastly, it is an attractive prospect to couple nonlinear waves to larger-scale circulation patterns \cite{pedlosky2013geophysical, dijkstra2005low, moore2024large, whiteheadenergy}, density-stratified environments \cite{sutherland2019recent, chini2022exploiting, camassa2012stratified, mac2022morphological}, or porous environments \cite{han2020dynamic, chiu2020viscous, mccurdy2022predicting, moore2023fluid} in order to assess wave feedback in these systems.

\section*{Acknowledgments}
N.J.~Moore recognizes support from the Office of Naval Research, award number N000142412617. B.~Foerster recognizes support from the Colgate University student-initiated summer fellowship.

\bibliography{bib}

\appendix
\section{Validation with two-mode exact solutions}\label{appA}
In the case $\nmodes =2$, it is possible to solve exactly for the Hamiltonian components $\Ham_2$ and $\Ham_3$, offering a way to validate the numerical computation of these quantities. In particular, consider arbitrary values of two wave modes $\uhat_1$ and $\uhat_2$ written in polar form:
\begin{align}
\label{bmark_uhat}
\uhat_1 = R_1 e^{i \theta_1} , \quad \uhat_2 = R_2 e^{i \theta_2}.
\end{align}
The energy constraint, \cref{Econd}, gives
\begin{equation}
2\pi (\abs{\uhat_1}^2 + \abs{\uhat_2}^2) = 2\pi (R_1^2 + R_2^2)= \En_0.
\end{equation}
Therefore, the moduli can then be written as
\begin{equation}
(R_1, R_2) = \sqrt{ \frac{\En_0}{2\pi} } \,\, (\cos(\phi),\sin(\phi)) \, ,
\end{equation}
for some $\phi \in [0,2\pi)$. 
Thus, to select a set of arbitrary modes $\uhat_1$ and $\uhat_2$ satisfying the energy constraint, it suffices to specify the values $\theta_1, \theta_2, \phi \in [0,2\pi)$. For example, these three angles could be selected randomly from a uniform distribution.
Calculating $\Ham_2$ from \cref{H2Pars} and $\Ham_3$ from \cref{H3_double} gives the exact solutions:
\begin{align}
\Ham_2  = 2 \pi \left( R_1^2 + 4 R_2^2 \right) \, , \qquad
\Ham_3 = 2\pi R_2 R_1^2 \cos(2\theta_1-\theta_2) \, .
\end{align}
These exact solutions can be used to validate the numerical computation of $\Ham_2$ and $\Ham_3$ for any values of $\theta_1, \theta_2, \phi \in [0,2\pi)$.

\end{document}